\documentclass{amsart}
\newtheorem{THM}{Theorem}[section]
\newtheorem{LMA}[THM]{Lemma}
\newtheorem{PROP}[THM]{Proposition}
\newtheorem{CORO}[THM]{Corollary}

\newtheorem{EG}[THM]{Example}

\numberwithin{equation}{section}

\usepackage{amsmath}
\usepackage{amssymb}
\usepackage{euscript}

\newcommand{\showon}{\begin{eqnarray}}
\newcommand{\showoff}{\end{eqnarray}}

\newcommand{\CC}{\mathbf{C}}
\newcommand{\M}{\EuScript{M}}
\newcommand{\A}{\EuScript{A}}
\newcommand{\U}{\EuScript{U}}
\newcommand{\HH}{\EuScript{H}}
\renewcommand{\L}{\EuScript{L}}
\newcommand{\G}{\mathbf{G}}
\newcommand{\ssll}{\mathfrak{sl}}

\newcommand{\x}{\mathbf{x}}
\newcommand{\z}{\mathbf{z}}
\newcommand{\drop}{\smallsetminus}
\newcommand{\spn}{\mathrm{span}}
\newcommand{\goesto}{\rightarrow}
\newcommand{\none}{\varnothing}

\begin{document}

\title[Represented matroid algebras]{Algebras related to matroids
represented in characteristic zero}
\author{David G. Wagner}
\address{Department of Combinatorics and Optimization\\
University of Waterloo\\
Waterloo, Ontario, Canada\ \ N2L 3G1}
\email{\texttt{dgwagner@math.uwaterloo.ca}}
\thanks{Research supported by the Natural
Sciences and Engineering Research Council of Canada under
operating grant OGP0105392.}

\begin{abstract}
Let $k$ be a field of characteristic zero.  We consider graded subalgebras
$A$ of $k[x_{1},\ldots,x_{m}]/(x_{1}^{2},\ldots,x_{m}^{2})$ generated
by $d$ linearly independent linear forms.   Representations of matroids
over $k$ provide a natural description of the structure of these algebras.
In return, the numerical properties of the Hilbert function of $A$
yield some information about the Tutte polynomial of the corresponding
matroid.  Isomorphism classes of these algebras correspond 
to equivalence classes of hyperplane arrangements under the action
of the general linear group.
\end{abstract}
\maketitle

\begin{center}{\textit{Dedicated to the memory of Fran\c{c}ois Jaeger.}}
\end{center}\vspace{5mm}

\section{Introduction}

We consider the following class of graded
algebras over a field $k$ of characteristic zero.  Let $B:=
k[x_{1},\ldots,x_{m}]/(x_{1}^{2},\ldots,x_{m}^{2})$
with the standard grading (so $B=\bigoplus_{j=0}^{m}B_{j}$ and
$\dim_{k}B_{j}=\binom{m}{j}$), and let $A=\bigoplus_{j=0}^{m}A_{j}$
be a subalgebra of $B$ generated by $d$ linearly independent
forms of degree one.  Two examples motivate the investigation
of such algebras.

\begin{EG}\textup{
Let $G$ be a finite undirected graph with $m$ edges, and orient each
edge arbitrarily.  Fixing a bijection between the edges of $G$
and the indeterminates $\{x_{j}\}$, we regard a linear form in
$B_{1}$ as a linear combination of the edges of $G$.  Let $A_{1}$
be the ``cycle-space'' of $G$  (that is, the subspace of $B_{1}$
consisting of linear combinations of the oriented edges satisfying
Kirchhoff's First Law:  at every vertex the net flux is zero),
and let $A$ be the subalgebra of $B$ generated by $A_{1}$. In
\cite{W1} it is shown that this construction may be symmetrized
to obtain a graded algebra $\Phi_{\cdot}(G,k)$ which is independent
of the choice of orientation of the edges of $G$, and
which is covariantly functorial with respect to graph morphisms.
Formally, $\Phi_{\cdot}(G,k)$ resembles a cohomology ring for the
graph $G$ with coefficients in the field $k$.}\end{EG}

\begin{EG}\textup{
Let $G$ be a connected complex semisimple Lie group, with
Borel subgroup $B$ and root system $\Delta$, and consider the
homogeneous manifold $X=G/B$  (the ``flag manifold'' of type $G$).
Postnikov, Shapiro, and Shapiro \cite{PSS} (see also Shapiro, Shapiro,
and Vainshtein \cite{SSV})
identify differential two-forms $\{\phi_{\alpha}:\ \alpha\in\Delta\}$
on $X$ such that $\phi_{-\alpha}=-\phi_{\alpha}$, $\phi_{\alpha}^{2}=0$,
and the $\phi_{\alpha}$ pairwise commute.  Any weight $\lambda$ of $G$
determines a holomorphic hermitian line bundle $L_{\lambda}$ on
$X$, and the curvature form $\Theta(L_{\lambda})$ of this line
bundle is a linear combination of the $\{\phi_{\alpha}:\ \alpha\in
\Delta\}$. The subalgebra $C(X)$ of the algebra of differential
forms on $X$ generated by the curvature forms $\Theta(L_{\lambda})$
is of the kind considered here, and the cohomology ring
$H^{\cdot}(X,\CC)$ is a quotient of $C(X)$.
}\end{EG}

In the next section we show that an isomorphism class of algebras
$A$ as above corresponds to a linear equivalence class of representations
of a matroid over the field $k$.  Equivalently, this corresponds to
an equivalence class of hyperlane arrangements $\HH\subset k^{d}$
under the action of the general linear group $GL(k^{d})$.  One
direction of this correspondence is immediate (Lemma $1.2$) while the
other requires substantial preliminaries (Theorem $1.10$).
We establish a deletion/contraction short exact sequence
which proves to be useful (Theorem $1.5$).  The Poincar\'e polynomial
of $A$ is a specialization of the Tutte polynomial of the corresponding
matroid.  
We present $A$ as a quotient of a polynomial ring modulo an explicitly
given ideal (Theorem $1.8$), and prove an analogue of half of the Strong
Lefschetz Theorem for these algebras (Theorem $1.12$).  In Section 2 we
discuss inequalities on the Hilbert function of $A$ derived from the
algebraic structure of $A$.   Having computed a few hundred random
examples, it seems that the Hilbert function of $A$ is logarithmically
concave, and we prove this generically and in the case $d=2$.  These
results go some way towards addressing Problems 6.8 and 6.10 of \cite{W1}.

\section{Algebraic Structure}

For a natural number $n$ we use the notation $[n]:=\{1,2,\ldots,n\}$.
For $0\leq j\leq m$, let $\Delta_{j}$ be the set of square-free
monomials $\x^{\alpha}$ of degree $j$ in $\{x_{1},\ldots,x_{m}\}$, so 
$\Delta:=\bigcup_{j=0}^{m}\Delta_{j}$ is a $k$-basis for $B$.
Endomorphisms of $B_{j}$ are represented by square matrices with rows
and columns indexed by $\Delta_{j}$.  A \emph{monomial matrix} has
exactly one nonzero entry in each row and each column.

\begin{LMA} The $k$-algebra automorphisms of $B$ form a group
$\mathrm{Aut}_{k}(B)$ which is isomorphic to the
group of monomial matrices acting on $B_{1}$
with respect to the basis $\Delta_{1}$.
\end{LMA}
\begin{proof}
Notice that if $f\in B_{1}$ is such that $f^{2}=0$ then
$f=cx_{j}$ for some $c\in k$ and $j\in[m]$.  Thus,
for any automorphism $\phi:B\goesto B$ there is a permutation
$\sigma:[m]\goesto[m]$ and nonzero scalars $c_{j}\in k$
such that $\phi(x_{j})=c_{j}x_{\sigma(j)}$ for all $j\in[m]$.
Conversely, any such choice of $\sigma$
and $\{c_{j}\}$ determines an automorphism of $B$.
\end{proof}

Let $M=(m_{ij})$ be a $d$-by-$m$ matrix over
$k$ for which the rowspace of $M$ is $A_{1}$.  (Henceforth we
identify row vectors of
length $m$ with linear combinations of the indeterminates $\{x_{j}\}$.)
Since $M$ determines $A$ we will often use the notation $A(M)$.
The linearly independent sets of columns of $M$ form the
independent sets of a matroid $\M$, and $M$ is a 
\textit{representation of $\M$ over $k$}.  (For background information 
on matroids consult Oxley \cite{Ox} or Welsh \cite{We}.)
Two representations $M$ and $N$ of $\M$ are \textit{linearly equivalent}
if there is a monomial matrix $P$ and an invertible matrix $Q$ such
that $QMP=N$.

\begin{LMA}  Let $M$ and $N$ be two $d$-by-$m$ matrices of rank $d$
over the field $k$.  If $M$ and $N$ are linearly equivalent
representations of the same matroid, then $A(M)\simeq A(N)$.
\end{LMA}
\begin{proof}
If $QMP=N$ with $Q$ invertible and $P$ a monomial matrix, then by
Lemma 1.1, $P$ determines a $k$-algebra automorphism of $B$ such that
$A_{1}(M)\simeq A_{1}(MP)=A_{1}(N)$.  Since $A(M)$ and $A(N)$ are
generated by linear forms, it follows that $A(M)$ and $A(N)$ are isomorphic
$k$-algebras.
\end{proof}
The converse of Lemma $1.2$ also holds (Theorem $1.10$), but the proof
relies on a presentation of $A(M)$ as a quotient of a polynomial
ring (Theorem $1.8$) which takes some work to derive.

Lemma $1.2$ has an interesting geometric interpretation;  see Orlik and
Terao \cite{OT} for background on hyperplane arrangements.

\begin{EG}\textup{
Let $\HH$ be a (nonreduced, central, essential) arrangement of $m$
hyperplanes in a $d$-dimensional $k$-vectorspace $V$.  Choose an
arbitrary basis $\EuScript{B}$ of $V^{*}$, an arbitrary enumeration
$\HH=\{H_{1},\ldots,H_{m}\}$ of $\HH$, and arbitrary linear forms
$\ell_{1},\ldots,\ell_{m}$ in $V^{*}$ such that $H_{j}=\ker(\ell_{j})$
for $j\in[m]$.  Writing each $\ell_{j}$ as a column vector with 
respect to the basis $\EuScript{B}$ determines a $d$-by-$m$ matrix $M$
of rank $d$.  If $N$ is another such matrix obtained from $\HH$ by 
different choices of basis, enumeration, and linear forms, then there
is an invertible $d$-by-$d$ matrix $Q$ (for the change of basis) and
an $m$-by-$m$ monomial matrix $P$ (for change of enumeration and
rescaling of linear forms) such that $QMP=N$.  Therefore, by Lemma 
$1.2$, the algebra $A(\HH):=A(M)$ is a well-defined invariant of the
hyperplane arrangement.  Moreover, if $\HH'$ is a hyperplane arrangement
which is equivalent to $\HH$ under the action of $GL(V)$, then the
corresponding matrices $M$ and $M'$ are linearly equivalent 
representations of the same matroid, and so $A(\HH')\simeq A(\HH)$.
}\end{EG}

Lemma 1.4 prepares for Theorem 1.5.

\begin{LMA}  Consider linear forms 
$f_{i}=x_{i}+\sum_{j=d+1}^{m}c_{ij}x_{j}$ in $B_{1}$ for $i\in[d]$, and
a polynomial $p(z_{1},\ldots,z_{d})$ in $k[z_{1},\ldots,z_{d}]$.
If $f_{1}p(f_{1},\ldots,f_{d})=\sum_{\alpha}s_{\alpha}\x^{\alpha}\neq 0$
then there is some $\x^{\alpha}\in\Delta$ which is divisible by $x_{1}$ and
such that $s_{\alpha}\neq 0$.
\end{LMA}
\begin{proof}
Since $f_{1}p(f_{1},\ldots,f_{d})\neq 0$, there is some 
$\x^{\beta}\in\Delta$ with $s_{\beta}\neq 0$.  Let $T$ be the set of
$j\in[m]$ such that $x_{j}$ divides $\x^{\beta}$, $c_{1j}\neq 0$, and 
the coefficient $w_{j}$ of $\x^{\beta}x_{j}^{-1}$ in 
$p(f_{1},\ldots,f_{d})$ is nonzero.  Thus, $s_{\beta}=\sum_{j\in T}
c_{1j}w_{j}$.  If $x_{1}$ divides $\x^{\beta}$ then the result is 
proved, so we may assume that $x_{1}$ does not divide $\x^{\beta}$, 
and hence that $1\not\in T$.  Since $T$ is not empty there is some $j\in T$;
now consider the monomial $\x^{\alpha}:=x_{1}\x^{\beta}x_{j}^{-1}$.  We claim 
that this occurs in $f_{1}p(f_{1},\ldots,f_{d})$ with coefficient
$s_{\alpha}=w_{j}$, which is nonzero.  But this is clear, since
in  $f_{1}p(f_{1},\ldots,f_{d})=\sum_{a=1}^{b}q_{a}(f_{2},\ldots,f_{d})
f_{1}^{a}$ the terms contributing to $s_{\alpha}\x^{\alpha}$ correspond
bijectively with the terms contributing to $s_{\beta}\x^{\beta}$ which
choose $x_{j}$ from some factor $f_{1}$.  The correspondence is made
simply by replacing $x_{1}$ by $x_{j}$ in each such term, and the ratio
of the coefficients of corresponding terms is $1:c_{1j}$.  
\end{proof}

For a $d$-by-$m$ matrix $M$ and $j\in[m]$, let $M\drop j$ be
the $d$-by-$(m-1)$ matrix obtained by deleting the $j$-th column
from $M$.  If this column is identically zero then $A(M\drop j)\simeq
A(M)$, as is easily seen.  As a result, we are free to assume that
$M$ has no zero columns in what follows.  If column $j$ of $M$ is
not zero then let $i\in[d]$ be the greatest index such that
$m_{ij}\neq 0$, and produce $M'$ by adding $-m_{ih}/m_{ij}$
times column $j$ to column $h$ of $M$, for each $h\in[m]$.
Finally, $M/j$ is the $(d-1)$-by-$(m-1)$ matrix obtained by
deleting the $i$-th row and $j$-th column from $M'$.

Theorem 1.5 is an analogue of the sequence (3.1) of \cite{W1}.
(The notation $A(M\drop j)(-1)$ merely indicates that the grading of
$A(M\drop 1)$ has been shifted up by one degree.)

\begin{THM}  Let $M$ be a $d$-by-$m$ matrix of rank $d$ over the
field $k$. For each $j\in[m]$ such that column $j$ of $M$ is not zero,
there is a short exact sequence of graded $k$-spaces
$$0\longrightarrow A(M\drop j)(-1)\longrightarrow
A(M)\stackrel{\pi}{\longrightarrow}A(M/j)\longrightarrow 0$$
in which $\pi$ is a $k$-algebra homomorphism.
\end{THM}
\begin{proof}
Replacing $M$, if necessary, by a linearly equivalent representation
of the same matroid, we may assume that $j=1$ and that $M$ has
the block structure $M=[I\ N]$ in which $I$ is the $d$-by-$d$
identity matrix.  Let $f_{1},\ldots,f_{d}$ be the rows of
$M$, let $f'_{1},\ldots,f'_{d}$ be the rows of $M\drop 1$, and
let $f''_{2},\ldots,f''_{d}$ be the rows of $M/1$.  There is certainly
an exact sequence
$$0\longrightarrow (f_{1})\longrightarrow
A(M)\longrightarrow A(M)/(f_{1})\longrightarrow 0$$
for the principal ideal $(f_{1})$ of $A(M)$.  It remains only to
establish isomorphisms $A(M/1)\simeq A(M)/(f_{1})$ and $
A(M\drop 1)(-1)\simeq(f_{1})$.

Now, since column $1$ of $M$ is zero except in row $1$, $f''_{i}=f_{i}$
for $2\leq i\leq d$;  thus,
there is a well-defined $k$-algebra homomorphism from $A(M)$ to $A(M/1)$
given by $f_{1}\mapsto 0$ and $f_{i}\mapsto f''_{i}$ for
$2\leq j\leq d$.  Clearly this is surjective and has kernel $(f_{1})$.
For the other isomorphism, notice that $f_{1}=x_{1}+f'_{1}$ and 
$f_{i}=f'_{i}$ for $2\leq i\leq d$.  Thus, 
$f_{1}^{a}=(f'_{1})^{a}+ax_{1}(f'_{1})^{a-1}$ for every natural number 
$a$; it follows that for any polynomial $p(z_{1},\ldots,z_{d})$,
$$p(f_{1},\ldots,f_{d})=p(f'_{1},\ldots,f'_{d})+x_{1}p'(f'_{1},\ldots,f'_{d}),$$
in which $p'(\z):=(\partial/\partial z_{1})p(\z)$.  Thus, the rule
$p(f_{1},\ldots,f_{d})\mapsto p'(f'_{1},\ldots,f'_{d})$ gives a 
well-defined $k$-linear homomorphism $\phi:A(M)\goesto A(M\drop 
1)(-1)$;  this is just the extraction of the coefficient of $x_{1}$ 
from $p(f_{1},\ldots,f_{d})$.  Since for every polynomial 
$q(\z)$ there is a polynomial $p(\z)$ such that
$(\partial/\partial z_{1})z_{1}p(\z)=q(\z)$, it follows that the
restriction of $\phi$ to $(f_{1})$ is surjective onto $A(M\drop 1)(-1)$.
Finally, Lemma $1.4$ shows that the restriction of $\phi$ to
$(f_{1})$ is injective, establishing the isomorphism $A(M\drop 1)(-1)
\simeq(f_{1})$.
\end{proof}

The \textit{Poincar\'e polynomial} of a finite-dimensional
graded $k$-space $A=\bigoplus_{j=0}^{m} A_{j}$ is
$P(A;t):=\sum_{j=0}^{m}(\dim_{k} A_{j})t^{j}$.  The coefficients of
$P(A;t)$ form the \emph{Hilbert function} of $A$.  The \emph{Tutte
polynomial} $T_{\M}(x,y)$ of a matroid $\M$ is the class of $\M$ in
the Grothendieck ring of the category of matroids;  see \cite{BO,Ox,We}.

\begin{CORO} Let $M$ be a $d$-by-$m$ matrix of rank $d$ with no zero
columns, representing the matroid $\M$ over the field $k$.
Then the Poincar\'e polynomial of $A(M)$ depends only on $\M$ and is
$P(\M;t):=P(A(M);t)=t^{m-d}T_{\M}(1+t,t^{-1})$.
\end{CORO}
\begin{proof}
From Theorem 1.5 we obtain the recursion
\begin{eqnarray}
P(A(M);t)=tP(A(M\drop j);t)+P(A(M/j);t)
\end{eqnarray}
for the Poincar\'e polynomials, with initial conditions
$P(A(M);t)=1$ if $d=0$, and $P(A(M);t)=1+t+\cdots+t^{m}$
if $d=1$ and $M$ has $m$ nonzero columns.  Defining
$\widetilde{P}(A(M);t):=t^{d-m}P(A(M);t)$ we have
$\widetilde{P}(A(M);t)=\widetilde{P}(A(M\drop j);t)+
\widetilde{P}(A(M/j);t)$.  When $d=0$ and $m=1$ we have
$\widetilde{P}(A(M);t)=t^{-1}$, and when $d=1$ and $m=1$ we have
$\widetilde{P}(A(M);t)=1+t$.  Since $T_{\M}(x,y)$ is the universal
Tutte-Grothendieck invariant of the category of matroids (see Brylawski
and Oxley \cite{BO}), it follows by induction on $d$ and $m$ that
$\widetilde{P}(A(M);t)=T_{\M}(1+t,t^{-1})$.
\end{proof}

For the next result we need some operations on sets of
columns of the matrix $M$.  For $S\subseteq[m]$, let $\spn_{k}(S)$
be the $k$-space spanned by the columns of $M$ in $S$, and let
$\overline{S}$ be the set of columns of $M$ contained in $\spn_{k}(S)$.
The \emph{rank} of $S$ is $r(S):=\dim_{k}\spn_{k}(S)$.
Let $I(S)$ be the lexicographically earliest basis of $\spn_{k}(S)$
contained in $\overline{S}$.  A column $j$ is \emph{externally active}
for $S$ if and only if $j\in\overline{S}\drop S$ and
$I(S\cup\{j\})=I(S)$;  Let $EA(S)$ be the set of columns externally
active for $S$, and let $\mathrm{ea}(S)$ be the cardinality of this set.

What follows is a new proof of Theorem 2 of Postnikov, Shapiro, and Shapiro
\cite{PSS}.

\begin{THM} Let $M$ be a $d$-by-$m$ matrix of rank $d$
with no zero columns, representing the matroid $\M$ over the field $k$.
For $0\leq j\leq m$, 
$\dim_{k} A_{j}(M)$ is the number of independent sets of $\M$
such that $m-\#S-\mathrm{ea}(S)=j$.
\end{THM}
\begin{proof}
The rank-polynomial expansion (see (6.12) of Brylawski and Oxley \cite{BO})
of $T_{\M}(x,y)$ is
$$T_{\M}(x,y)=\sum_{S\subseteq [m]}(x-1)^{d-r(S)}(y-1)^{\#S-r(S)}.$$
Making the substitution of Corollary $1.6$ leads to
\begin{eqnarray*}
P(\M;t) &=& \sum_{S\subseteq[m]}t^{m-\#S}(1-t)^{\#S-r(S)}
= \sum_{S\subseteq[m]}t^{m-\#S}\sum_{T\subseteq S\drop I(S)}
(-1)^{\#T}t^{\#T}\\
&=& \sum_{R\subseteq[m]}t^{m-\#R}\sum_{T\subseteq EA(R)}(-1)^{\#T}
=\sum_{R\subseteq[m]:\ EA(R)=\none}t^{m-\#R}.
\end{eqnarray*}
For the third equality, notice that $S\drop I(S)\subseteq EA(I(S))$ 
for every $S\subseteq[m]$.  Thus, if $T\subseteq S\drop I(S)$ and 
$R:=S\drop T$ then $I(R)=I(S)$ and $T\subseteq EA(R)$.  Conversely,
if $T\subseteq EA(R)$ then $T\subseteq (R\cup T)\drop I(R\cup T)$.

Since $EA(R)=EA(I(R))\drop R$, it follows that $EA(R)=\none$ if and only if
$R=I(R)\cup EA(I(R))$.  Conversely, if $S$ is independent then
$I(S\cup EA(S))=S$. Thus, the functions $R\mapsto I(R)$ and $S\mapsto
S\cup EA(S)$ are mutually inverse bijections between the sets
$\{R\subseteq[m]:\ EA(R)=\none\}$ and $\{S\subseteq[m]:\ S\ \mathrm{is\
independent\ in}\ \M\}$.  Therefore,
$$P(A(M);t) = \sum_{R\subseteq[m]:\ EA(R)=\none}t^{m-\#R}
= \sum_{S} t^{m-\#S-\mathrm{ea}(S)},$$
with the last sum over the independent sets of $\M$.
\end{proof}

We next present $A(M)$ as a quotient of the polynomial ring
$R:=k[z_{1},\ldots,z_{d}]$.  For any linear form $f=\sum_{j=1}^{m}
c_{j}x_{j}$ in $B_{1}$, let $\nu(f):=\#\{j:\ c_{j}\neq 0\}$.  Notice
that $f^{\nu(f)}\neq 0$ and $f^{1+\nu(f)}=0$.  Identifying a linear
form $p$ in $R_{1}$ with a row vector of length $d$, there is a
corresponding linear form $pM$ in $B_{1}$.  Define the ideal $J(M)$ of
$R$ by
$$J(M):=(p^{1+\nu(pM)}:\ p\in R_{1}).$$
Theorem 1.8 generalizes Theorem 4.8 of \cite{W1} and Proposition 1.1 of 
Shapiro, Shapiro, and Vainshtein \cite{SSV}.

\begin{THM} For $M$ a $d$-by-$m$ matrix of rank $d$ over the field
$k$, $A(M)\simeq R/J(M)$.
\end{THM}
\begin{proof}
We apply Theorem 1.5 for some $j\in[m]$
indexing a nonzero column of $M$.  Replacing $M$, if necessary, by a
linearly equivalent representation of the same matroid, we may assume
that $j=1$ and that $M$ has the block structure $M=[I\ N]$ in which $I$
is the $d$-by-$d$ identity matrix.  Let $f_{1},\ldots,f_{d}$ be the rows
of $M$, let $f'_{1},\ldots,f'_{d}$ be the rows of $M\drop 1$, and
let $f''_{2},\ldots,f''_{d}$ be the rows of $M/1$. 

Define a $k$-algebra homomorphism $\psi:R\goesto A(M)$ by
$\psi(z_{i}):=f_{i}$ for $i\in[d]$.  Certainly $\psi$
is surjective, as $A_{1}$ generates $A$.  We claim that
$\ker(\psi)=J(M)$, which we prove by induction on $d$ and $m$, the bases
$d=1$ and $m=d$ being easily seen.
It is clear that $J(M)\subseteq \ker(\psi)$ since for any
$p\in R_{1}$ we have $\psi(p)=pM$ and $(pM)^{1+\nu(pM)}=0$ in $A(M)$.
For the converse, define a $k$-algebra homomorphism
$\psi':R\goesto A(M\drop 1)$ by $\psi'(z_{i}):=f'_{i}$ for $i\in[d]$,
and define $\psi'':k[z_{2},\ldots,z_{d}]\goesto A(M/1)$ by
$\psi''(z_{i}):=f''_{i}$ for $2\leq i\leq d$.  There is a commutative diagram
$$\begin{array}{ccccccccc}
0 & \longrightarrow & R(-1) &
\stackrel{\eta}{\longrightarrow} & R & \stackrel{\pi}{\longrightarrow}
& k[z_{2},\ldots,z_{d}] & \longrightarrow & 0 \\
 & & \downarrow\scriptsize{\psi'} & & \downarrow\scriptsize{\psi} & &
 \downarrow\scriptsize{\psi''} & & \\
0 & \longrightarrow & A(M\drop 1)(-1) &\longrightarrow &
A(M) &\longrightarrow& A(M/1) &
\longrightarrow & 0
\end{array}$$
in which the bottom row is the sequence of Theorem 1.5.  From the proof
of Theorem $1.5$ one sees that the homomorphisms in the top row are
given by $\pi(p(\z)):=p(0,z_{2},\ldots,z_{d})$ and
$\eta(p(\z)):=\int p(\z)\mathrm{d}z_{1}$
for all $p(\z)\in R$.  Since $\psi'$ is surjective, 
the kernel-cokernel exact sequence (see, \textit{e.g.} Lemma II.5.2 of
Mac Lane \cite{Mac}) implies that
$0\goesto \ker(\psi')\goesto \ker(\psi)\goesto \ker(\psi'')\goesto 0$
is exact.  By induction, we deduce that
$$\ker(\psi)=\eta(J(M\drop 1)(-1))\oplus \iota(J(M/1)),$$
in which $\iota:k[z_{2},\ldots,z_{d}]\goesto R$ is the natural inclusion.

To prove that $\ker(\psi)\subseteq J(M)$, it thus suffices to show
that $\iota(J(M/1))\subseteq J(M)$ and $\eta(J(M\drop 1)(-1))\subseteq
J(M)$.  The first of these claims is trivial, since $\iota(J(M/1))$
consists of exactly those polynomials in $J(M)$ which do not involve
the indeterminate $z_{1}$.  For the second claim, by $k$-linearity it
suffices to prove that $\eta(\z^{\gamma}g(\z))\in J(M)$
for any monomial $\z^{\gamma}$ and generator $g(\z)$ of $J(M\drop 1)(-1)$.
So, let $p(\z):=c_{1}z_{1}+\cdots+c_{d}z_{d}$ and let $\nu:=
\nu(c_{1}f'_{1}+\cdots+c_{d}f'_{d})$, and consider 
$\eta(\z^{\gamma}p(\z)^{1+\nu})$.  
If $c_{1}=0$ then $\nu(c_{1}f_{1}+\cdots+c_{d}f_{d})=\nu$ and 
$\int \z^{\gamma}p(\z)^{1+\nu}\mathrm{d}z_{1}=
\z^{\gamma}z_{1}p(\z)^{1+\nu}/(\gamma_{1}+1)$
is in $J(M)$.  On the other hand, if $c_{1}\neq 0$ then
$\nu(c_{1}f_{1}+\cdots+c_{d}f_{d})=\nu+1$;  however, applying 
integration by parts repeatedly we obtain
\begin{eqnarray*}
\int \z^{\gamma}p(\z)^{1+\nu}\mathrm{d}z_{1} &=&
\frac{\z^{\gamma}p(\z)^{2+\nu}}{2+\nu}-\int\left(\frac{\partial\z^{\gamma}}
{\partial z_{1}}\right)\frac{p(\z)^{2+\nu}}{2+\nu}\mathrm{d}z_{1}\\
&=& \cdots\ = q(\z)p(\z)^{2+\nu}
\end{eqnarray*}
for some polynomial $q(\z)\in R$.  Since $p(\z)^{2+\nu}$ is a 
generator of $J(M)$, the result follows.
\end{proof}
Although Theorem 1.8 gives a good picture of $A(M)$, it 
would be preferable to have a standard monomial theory for this
algebra.  Presumably this would rely on matroid-theoretic structure
as in the proof of Theorem 1.7, but as yet the situation remains 
unclear.

We can now establish the converse of Lemma $1.2$, the proof of which uses
the following ``tomographic'' lemma (valid for any infinite
field $k$).

\begin{LMA}  Let $\L$ and $\L'$ be finite multisets of lines in the
$d$-dimensional $k$-vectorspace $V$, each line passing through the
origin.  Assume that for every hyperplane $H\subset V$, the number
of lines of $\L$ in $H$ equals the number of lines of $\L'$ in $H$.
Then $\L=\L'$.
\end{LMA}
\begin{proof}  Arguing by contradiction, assume that $\L\neq\L'$.
Replacing $\L$ by $\L\drop\L'$ and $\L'$ by $\L'\drop\L$, we may assume
that $\L\cap\L'=\none$.  At least one of $\L$ or $\L'$ is nonempty;
by symmetry, consider any $\ell\in\L$.  Since $k$ is infinite, there
are infinitely many hyperplanes $H\subset V$ containing $\ell$.  By
the hypothesis, each of these hyperplanes contains at least one line
from $\L'$, which is not $\ell$. Since these lines must be pairwise
distinct $\L'$ is infinite, a contradiction.
\end{proof}

\begin{THM}  Let $M$ and $N$ be two $d$-by-$m$ matrices of rank $d$
over the field $k$.  Then $A(M)$ and $A(N)$ are isomorphic as
$k$-algebras if and only if $M$ and $N$ are linearly equivalent
representations of the same matroid.
\end{THM}
\begin{proof}
Lemma $1.2$ establishes one direction.  For the converse, assume
that $\phi:A(M)\goesto A(N)$ is a $k$-algebra isomorphism.
By the remarks preceding Theorem $1.5$, we may assume that $M$ and
$N$ have no zero columns.  Let $f_{1},\ldots,f_{d}$ be the rows of
$M$, and let $g_{1},\ldots,g_{d}$ be the rows of $N$.  Replacing $N$,
if necessary, by a linearly equivalent representation of the same
matroid, we may assume that $\phi:A_{1}(M)\goesto A_{1}(N)$ is determined
by $\phi(f_{i})=g_{i}$ for all $i\in[d]$.  Now let
$R:=k[z_{1},\ldots,z_{d}]$ and define
$\psi:R\goesto A(M)$ and $\psi':R\goesto A(N)$ by $\psi(z_{i}):=f_{i}$
and $\psi'(z_{i}):=g_{i}$ for all $i\in[d]$.  From Theorem $1.8$ it
follows that $J(M)=\ker(\psi)=\ker(\psi')=J(N)$.
Let $\L$ be the multiset of lines in $k^{d}$ consisting of the scalar
multiples of the columns of $M$.  Let $\L'$ be the corresponding
multiset of lines for $N$.  Since $J(M)=J(N)$, for any linear form
$p\in R_{1}$ it follows that $\nu(pM)=\nu(pN)$;  that is, the
number of lines of $\L$ in $\ker(p)$ equals the number of lines of
$\L'$ in $\ker(p)$.  By Lemma $1.9$ it follows that $\L=\L'$.  Thus,
there is an $m$-by-$m$ monomial matrix $P$ such that $MP=N$.
This completes the proof.
\end{proof}

\begin{CORO}
Let $A$ and $A'$ be subalgebras of $B$ generated by linear forms.
Any $k$-algebra isomorphism $\phi:A\goesto A'$ extends
to an automorphism of $B$.
\end{CORO}
\begin{proof}
Let $A=A(M)$ and $A'=A(N)$ for $d$-by-$m$ matrices $M$ and $N$ of
rank $d$.  By Theorem $1.10$ there is an $m$-by-$m$ monomial
matrix $P$ such that $A_{1}(MP)=A_{1}(N)$ and for $f\in A_{1}(M)$,
$\phi(f)=fP$.  By Lemma $1.1$, $P$ determines an automorphism of
$B$ extending $\phi:A\goesto A'$.
\end{proof}

We close this section with an analogue of half of the Strong
Lefschetz Theorem, generalizing Theorem 4.10 of \cite{W1}.

\begin{THM} Let $M$ be a $d$-by-$m$ matrix of rank $d$
over the field $k$, and assume that $M$ has no zero columns.
Let $g=\sum_{j=1}^{m}c_{j}x_{j}\in A_{1}(M)$ be such that
$c_{j}\neq 0$ for all $j\in[m]$.
Then for each $0\leq j\leq m/2$, the homomorphism
$\cdot g^{m-2j}:A_{j}(M)\goesto A_{m-j}(M)$ is injective.
\end{THM}
\begin{proof}
 Fix any
$0\leq j\leq m/2$.  Let $W$ be the matrix with rows indexed
by $\Delta_{m-j}$ and with columns indexed by $\Delta_{j}$,
with $W_{\x^{\alpha},\x^{\beta}}:=1$ if $\x^{\beta}$ divides $\x^{\alpha}$
and $W_{\x^{\alpha},\x^{\beta}}:=0$ otherwise.  Multiplication by the element
$g^{m-2j}$ of $A(M)$ induces a homomorphism
$\cdot g^{m-2j}:B_{j}(M)\goesto B_{m-j}(M)$;  let $G$ be the
matrix representing this homomorphism with repsect to the bases
$\Delta_{j}$ and $\Delta_{m-j}$.  That is, if $\x^{\beta}|\x^{\alpha}$
then $G_{\x^{\alpha},\x^{\beta}}=\prod\{c_{j}:\ x_{j}|\x^{\alpha-\beta}\}$,
and $G_{\x^{\alpha},\x^{\beta}}:=0$ otherwise.  Let $P$ be the diagonal
square matrix indexed by $\Delta_{m-j}$, with $P_{\x^{\beta},\x^{\beta}}:=
\prod\{c_{j}:\ x_{j}|\x^{\beta}\}$, and let $Q$ be the diagonal
square matrix indexed by $\Delta_{j}$, with $Q_{\x^{\alpha},\x^{\alpha}}:=
\prod\{c_{j}:\ x_{j}|\x^{\alpha}\}$.  By the hypothesis on $g$,
both $P$ and $Q$ are invertible.  One verifies that the matrix
equation $GP=QW$ holds, and hence $\det(G)=\det(Q)\det(W)
\det(P)^{-1}$.  Wilson \cite{Wil} proves that
$$\det(W)=\prod_{h=0}^{j}\binom{m-j-h}{j-h}^{\binom{m}{h}-
\binom{m}{h-1}},$$
and therefore $\det(G)\neq 0$.  Thus,
$\cdot g^{m-2j}:B_{j}(M)\goesto B_{m-j}(M)$ is an isomorphism,
and so $\cdot g^{m-2j}:A_{j}(M)\goesto A_{m-j}(M)$ is injective.
\end{proof}

\section{Numerics of the Hilbert functions}

The notation $d_{j}(M):=\dim_{k}A_{j}(M)$ is convenient for this
section. 
For positive integers $a$ and $j$ there is a unique expression 
$$a=\binom{a_j}{j}+\binom{a_{j-1}}{j-1}+\cdots+\binom{a_i}{i}$$
such that $a_j>a_{j-1}>\cdots>a_i\geq i>0$.  The
\emph{$j$-th pseudopower of $a$} is 
$$\psi_j(a):=\binom{a_j+1}{j+1}+\binom{a_{j-1}+1}{j}+\cdots+
\binom{a_i+1}{i+1}.$$

\begin{CORO}
Let $M$ be a $d$-by-$m$ matrix of rank $d$ with no zero columns.
Then $d_{0}(M)=1$, $d_{1}(M)=d$, $d_{m}(M)=1$, and for $j\in[m-1]$,
we have $0<d_{j+1}(M)\leq\psi_j(d_j(M))$.
\end{CORO}
\begin{proof}
Since $A_{0}(M)=k$ and $\dim_{k}A_{1}(M)=d$, the first two statements
are clear.  Since $M$ has no zero columns and $k$ is infinite, there is
a linear form $g\in A_{1}$ with $\nu(g)=m$.  Thus, $g^{m}$ is a
nonzero multiple of $x_{1}\cdots x_{m}$; since $B_{m}$ is
one-dimensional it follows that $d_{m}(M)=1$.  The remaining 
inequalities are a direct application of Macaulay's Theorem (see Theorems
II.2.2 and II.2.3 of Stanley \cite{St}), since $A$ is generated by 
linear forms.
\end{proof}

\begin{CORO}
Let $M$ be a $d$-by-$m$ matrix of rank $d$ with no zero columns.
Then $d_0(M)\leq d_1(M)\leq\cdots\leq d_{\lfloor m/2\rfloor}(M)$, and if
$0\leq j\leq m/2$ then $d_j(M)\leq d_{m-j}(M)$.
\end{CORO}
\begin{proof} Since $M$ has no zero columns and $k$ is infinite, there is
a linear form $g\in A_{1}$ with $\nu(g)=m$.  The monomorphisms
$\cdot g^{m-2j}:A_{j}(M)\goesto A_{m-j}(M)$ of Theorem $1.12$ show that
$d_{j}(M)\leq d_{m-j}(M)$ for all $0\leq j\leq m/2$.  Since each of
these maps is injective, each of the maps $\cdot g:A_{j}(M)\goesto 
A_{j+1}(M)$ for $0\leq j<m/2$ must also be injective, implying the 
remaining inequalities.
\end{proof}

Notice that
$d_{j}(M)\leq\binom{d+j-1}{j}$ for all
$j\geq 0$ since $A$ is generated by $d$ linear forms, and that
$d_{j}(M)\leq\binom{m}{j}$ for all $j\geq 0$ since $A$ is a
subalgebra of $B$.  Next, we see that generically these bounds are 
attained.

\begin{PROP} Let $M$ be a $d$-by-$m$ matrix
over the field $k$ representing the uniform matroid $\U_{m}^{d}$
of rank $d$ on $m$ elements.  Then for $0\leq j\leq m$,
$d_{j}(M)=\min\{\binom{d+j-1}{j},\binom{m}{j}\}$.
\end{PROP}
\begin{proof}
When $d\geq 2$, $M\drop 1$ represents $\U_{m-1}^{d}$ and $M/1$ 
represents $\U_{m-1}^{d-1}$.  Thus, from $(1.1)$ we obtain
$$P(\U_{m}^{d};t)=tP(\U_{m-1}^{d};t)+P(\U_{m-1}^{d-1};t)$$
for $d\geq 2$, with initial conditions $P(\U_{m}^{1};t)=1+t+\cdots+t^{m}$.
The result follows by induction on $d$ and $m$, using familiar
recurrences for binomial coefficients.
\end{proof}

Let $\A(m,d)$ be the collection of all graded subalgebras $A\subseteq B$
generated by $A_{1}$ and with $\dim_{k} A_{1}=d$.
This $\A(m,d)$ is a sub-bundle of the
trivial vector-bundle $\G(B_{1},d)\times B$ over the Grassmann
variety $\G(B_{1},d)$ of $d$-dimensional subspaces of $B_{1}$;  for a given
$d$-plane $A_{1}\subseteq B_{1}$, the fibre over $A_{1}$ is the
subalgebra of $B$ generated by $A_{1}$.  As the Poincar\'e
polynomial $P(A;t)$ varies with $A_{1}$, the rank of $\A(m,d)$ is not 
constant, so $\A(m,d)$ is not complete.
By upper semicontinuity, the rank of the fibre of $\A(m,d)$ over
$A_{1}(M)$ attains its generic value if the matroid represented by
$M$ is uniform.  We next prove the converse, giving equations in
local coordinates for the degeneracy locus of $\A(m,d)$.  Consider the
affine open chart $\EuScript{C}\subset\G(B_{1},d)$ of
$d$-planes $A_{1}\subseteq B_{1}$ of the form $A_{1}(M)$ for a $d$-by-$m$
matrix $M=[I\ N]$ with $I$ the $d$-by-$d$ identity matrix.
The entries of $N=(n_{ij})$ are local coordinates 
on $\EuScript{C}$.  Since $\G(B_{1},d)$ is covered by affine opens which
are in the orbit of $\EuScript{C}$ under Aut$_{k}(B)$, it suffices to
consider just this one chart $\EuScript{C}$.

\begin{PROP} Let $M=[I\ N]$ be a $d$-by-$m$ matrix with $I$ the 
$d$-by-$d$ identity matrix, representing the matroid $\M$ over the 
field $k$. The following are equivalent:\\
\textup{(a)}\  $\M$ is not the uniform matroid $\U_{m}^{d}$.\\
\textup{(b)}\  For some $j\in[m]$, $d_{j}(M)<
\min\{\binom{d+j-1}{j},\binom{m}{j}\}$.\\
\textup{(c)}\  For some $h\in[\min\{d,m-d\}]$ and some $h$-by-$h$ 
submatrix $N'$ of $N$, $\det(N')=0$.
\end{PROP}
\begin{proof}
Proposition $2.3$ shows that (b) implies (a).  To see that (a)
implies (c), if $\M$ is not uniform then there is a $d$-by-$d$
submatrix $M'$ of $M$ which is singular.  Deleting the rows and
columns of $M'$ which contain a $1$ from the $I$ block of $M$
produces a singular square submatrix of $N$, proving (c).  This
argument may be reversed to show that (c) implies (a) as well.  Finally,
assuming (a), if $\M$ is not $\U_{m}^{d}$ then from Theorem 1.7,
with $S$ ranging over the independent sets of $\M$,
$$P(A(M);1)=\sum_{S}1<\binom{m}{0}+\binom{m}{1}+
\cdots+\binom{m}{d}=\sum_{j=0}^{m}\min\left\{\binom{d+j-1}{j},\binom{m}{j}
\right\},$$
and (b) follows.
\end{proof}

A sequence $(d_{0},\ldots,d_{m})$ of positive integers is
\emph{logarithmically concave} if $d_{j}^{2}\geq d_{j-1}d_{j+1}$ for
all $2\leq j\leq m-1$.  From Proposition $2.3$ it follows that if
$M$ represents the uniform matroid $\U_{m}^{d}$ then the Hilbert
function of $A(M)$ is logarithmically concave.  (The argument is easy:
each of the sequences $\binom{d+j-1}{j}$ and $\binom{m}{j}$ for $0\leq
j\leq m$ is logarithmically concave, and the coefficientwise minimum 
of two logarithmically concave sequences is also logarithmically
concave.)  Thus, generically, the sequence of ranks of the graded
pieces of $\A(m,d)$ is logarithmically concave.  Whether or not this
remains true over the degeneracy locus of $\A(m,d)$ is an interesting
question;  one possible approach is as follows.

As observed in \cite{W2}, the Hilbert function of a graded $\CC$-space
$A=\bigoplus_{j=0}^{m} A_{j}$ is logarithmically concave if and only
if there is a representation of $\ssll_{2}(\CC)$ on $A\otimes A$ for
which the standard basis elements $\{\mathsf{X}, \mathsf{Y}, \mathsf{H}\}$ of 
$\ssll_{2}(\CC)$ act such that
$\mathsf{X}:A_{i}\otimes A_{j}\rightarrow A_{i-1}\otimes A_{j+1}$ and
$\mathsf{Y}:A_{i}\otimes A_{j}\rightarrow A_{i+1}\otimes A_{j-1}$ for
all $i$ and $j$.  Hence, such a representation exists on the generic
fibre of $\A(m,d)\otimes\A(m,d)$.  The difficulty lies in
degenerating this generic representation over $\mathrm{Spec}\,\CC[u]$ so 
that at $u=0$ a representation on the fibre above an arbitrary point
of $\G(B_{1},d)$ is obtained.  It is not clear how (or whether!) this
can be done, but the following degenerations of the irreducible
representations of $\ssll_{2}(\CC)$ seem relevant.  For a proposition
$P$, let $\langle P\rangle$ be $1$ if $P$ is true and $0$ if $P$ is 
false.  For integers $1\leq r\leq n$ of the same parity, let 
$X_{n,r}(u)$ be the $n$-by-$n$ matrix with entries
$$X_{n,r}(u)_{ij}:=\left\{\begin{array}{ll}
iu^{\langle i\leq r\rangle-\langle j> n-r\rangle} & \mathrm{if}\ i=j-1,\\
0 & \mathrm{otherwise,}\end{array}\right.$$
let $Y_{n,r}(u)$ be the $n$-by-$n$ matrix with entries
$$Y_{n,r}(u)_{ij}:=\left\{\begin{array}{ll}
(n-j)u^{\langle i>n-r\rangle-\langle j\leq r\rangle} & \mathrm{if}\ i=j+1,\\
0 & \mathrm{otherwise,}\end{array}\right.$$
and let $H_{n,r}(u)=X_{n,r}(u)Y_{n,r}(u)-Y_{n,r}(u)X_{n,r}(u)$.
For example, with $n=5$ and $r=3$,
$$X_{5,3}(z):=\left[\begin{array}{ccccc}
0 & u & 0 & 0 & 0\\
0 & 0 & 2 & 0 & 0\\
0 & 0 & 0 & 3 & 0\\
0 & 0 & 0 & 0 & 4u^{-1}\\
0 & 0 & 0 & 0 & 0\end{array}\right]\ \ \ \mathrm{and}\ \ \
Y_{5,3}(z):=\left[\begin{array}{ccccc}
0 & 0 & 0 & 0 & 0\\
4u^{-1} & 0 & 0 & 0 & 0\\
0 & 3 & 0 & 0 & 0\\
0 & 0 & 2 & 0 & 0\\
0 & 0 & 0 & u & 0\end{array}\right].$$
For $0\neq u\in\CC$ these matrices define an irreducible representation
of $\ssll_{2}(\CC)$ on $\CC^{n}$.  As $u\goesto 0$ these linear
transormations cease to be defined on all of $\CC^{n}$.  At
$u=0$ they remain defined on an $r$-dimensional subspace, on which they
still provide an irreducible representation of $\ssll_{2}(\CC)$.

In the special case $d=2$ we can establish a property stronger than 
logarithmic concavity by other means.

\begin{THM}  Let $M$ be a $2$-by-$m$ matrix of rank $2$
over the field $k$.
Then the sequence $d_{j}(M)-d_{j-1}(M)$ is
nonincreasing as $j$ goes from $1$ to $m$.  Consequently, the
Hilbert function of $A(M)$ is logarithmically concave.
\end{THM}
\begin{proof}
We may assume that $M=[I\ N]$ with $I$ the $2$-by-$2$ identity
matrix, and denote the rows of $M$ by $f_{1}$ and $f_{2}$.  By
Theorem 1.8, $A(M)\simeq R/J(M)$ in which $R=k[z_{1},z_{2}]$ and 
$$J(M):=((c_{1}z_{1}+c_{2}z_{2})^{1+\nu(c_{1}f_{1}+c_{2}f_{2})}:\
c_{1},c_{2}\in k).$$
The columns of $M$ are partitioned uniquely into subsets $E_{1},\ldots,
E_{s}$ such that columns $j$ and $j'$ belong to the same part $E_{h}$
if and only if they are proportional; for each $h\in[s]$,
let $e_{h}:=\#E_{h}$.  For each $h\in[s]$ there is a particular ratio
$c_{1}:c_{2}$ such that the $j$-th entry of $c_{1}f_{1}+c_{2}f_{2}$
is zero if and only if $j\in E_{h}$.  Thus, for each $h\in[s]$ there is
a linear form $p_{h}\in R$ such that $\nu(p_{h}M)=m-e_{h}$, and in fact
$J(M)$ is generated by $\{p_{1}^{1+m-e_{1}},\ldots,p_{s}^{1+m-e_{s}}\}$.
For each $0\leq j\leq m$ let $w_{j}:=\#\{h\in[s]:\ 1+m-e_{h}=j\}$.

Now $\dim_{k} R_{j} = j+1$ for all $j\geq 0$, and $\dim_{k} A_{m}(M)=1$
by Corollary $2.1$, so $\dim_{k} J_{m}(M)=m$.  But
$$\dim_{k} J_{m}(M)\leq\sum_{h=1}^{s}\dim_{k}(p_{h}^{1+m-e_{h}})_{m}
=\sum_{h=1}^{s}\left[m-(1+m-e_{h})+1\right]=m,$$
and since equality holds the forms $p_{h}^{1+m-e_{h}}$ for $h\in[s]$
impose independent conditions on homogeneous $j$-forms in $R$, for all
$0\leq j\leq m$.  Thus, for each $0\leq j\leq m$,
$$\dim_{k}A_{j}(M)=j+1-\sum_{i=0}^{j}w_{i}(j-i+1).$$
From this the inequalities $d_{j}(M)-d_{j-1}(M)\geq 
d_{j+1}(M)-d_{j}(M)$ for $j\in[m-1]$ follow.  By the inequality of
arithmetic and geometric means it follows that
$d_{j}(M)\geq(d_{j-1}(M)+d_{j+1}(M))/2\geq(d_{j-1}(M)d_{j+1}(M))^{1/2}$,
completing the proof.
\end{proof}

\vspace{1cm}
\begin{center}\textsc{Acknowledgments}\end{center}

I thank Karen Chandler, Aldo Conca, James Geelen, and Alexander
Postnikov for helpful conversations on various aspects of this work.
I am especially grateful to Leslie Roberts for pointing out an error
in an earlier proof of Theorem 1.10.

\end{document}